\documentclass[12pt,twoside]{amsart}
\usepackage{amssymb}
\usepackage{amscd}
\usepackage{xypic}

\title[Base point free theorems]{Base point free theorems\\
---saturation, b-divisors, and canonical bundle formula---}
\author{Osamu Fujino} 
\date{2011/2/17, version 1.15}
\address{Department of Mathematics, Faculty of Science, 
Kyoto University, Kyoto 606-8502, Japan}
\email{fujino@math.kyoto-u.ac.jp}
\newcommand{\Nklt}[0]{{\operatorname{Nklt}}}
\newcommand{\Div}[0]{{\operatorname{Div}}}
\newcommand{\Supp}[0]{{\operatorname{Supp}}}
\newcommand{\Exc}[0]{{\operatorname{Exc}}}

\newcommand{\rank}[0]{{\operatorname{rank}}}
\newcommand{\red}[0]{{\operatorname{red}}}
\newcommand{\xprojlim}[0]{{\operatorname{projlim}}}
\newcommand{\Spec}[0]{{\operatorname{Spec}}}

\newtheorem{thm}{Theorem}[section]
\newtheorem{lem}[thm]{Lemma}
\newtheorem{cor}[thm]{Corollary}
\newtheorem{prop}[thm]{Proposition}

\theoremstyle{definition}
\newtheorem{defn}[thm]{Definition}
\newtheorem{ex}[thm]{Example}

\newtheorem{rem}[thm]{Remark}
\newtheorem*{ack}{Acknowledgments}       
 
\newtheorem*{notation}{Notation}

\newtheorem{say}[thm]{}
\begin{document}
\bibliographystyle{amsalpha+}

\begin{abstract}
We reformulate base point free theorems. 
Our formulation is flexible 
and has some important applications. 
One of the main purposes of this paper 
is to prove a generalization of 
the base point free theorem in Fukuda's paper:~On numerically 
effective log canonical 
divisors. 
\end{abstract}

\maketitle 
\tableofcontents
\section{Introduction} 
In this paper, we reformulate base point free theorems by using 
Shokurov's ideas:~{\em{b-divisors, saturation 
of linear systems}}. 
Combining the refined Kawamata--Shokurov base 
point free theorem (cf.~Theorem \ref{thm1}) or its generalization 
(cf.~Theorem \ref{bigfukuda}) with 
Ambro's formulation of Kodaira's canonical bundle formula, 
we obtain various new base point free theorems 
(cf.~Theorems \ref{thm4}, \ref{abundantfukuda}). 
They are flexible and have some important applications (cf.~Theorem \ref{quasi1}). 
One of the main purposes of this paper is 
to prove a generalization of the base point free theorem 
in Fukuda's paper \cite{fukuda2}:~On numerically 
effective log canonical 
divisors. See \cite[Proposition 3.3]{fukuda2}.  

\begin{thm}[{cf.~Corollary \ref{fuku-main}}]
Let $(X, B)$ be an lc pair and 
let $\pi:X \to S$ be a proper 
morphism onto a variety $S$. Assume the following 
conditions{\em{:}} 
\begin{itemize}
\item[(a)] $H$ is a $\pi$-nef $\mathbb Q$-Cartier $\mathbb Q$-divisor on $X$, 
\item[(b)] $H-(K_X+B)$ is $\pi$-nef and $\pi$-abundant, 
\item[(c)] $\kappa (X_\eta, (aH-(K_X+B))_\eta)\geq 0$ 
and $\nu(X_\eta, (aH-(K_X+B))_\eta)=\nu(X_\eta, (H-(K_X+B))_\eta)$ for 
some $a\in \mathbb Q$ with $a>1$, 
where $\eta$ is the generic point of $S$,  
\item[(f)] there is a positive integer $c$ such that 
$cH$ is Cartier and that $\mathcal O_T(cH):=\mathcal O_X(cH)|_T$ is 
$\pi$-generated, 
where $T=\Nklt(X, B)$ is the non-klt locus of $(X, B)$.  
\end{itemize}
Then $H$ is $\pi$-semi-ample. 
\end{thm} 

We need it to prove the finite generation of the log canonical 
ring for log canonical $4$-folds in \cite{fujino-finite}. 
See \cite[Remark 3.4]{fujino-finite}.  
As we explained in \cite[Remark 3.10.3]{fujino} and \cite[5.1]{fujino-kawamata}, 
the proof of \cite[Theorem 4.3]{kawamata} contains a gap. 
There are no rigorous proofs of \cite[Theorem 5.1]{kawamata} by the gap 
in the proof of  \cite[Theorem 4.3]{kawamata}, and 
the proof of \cite[Proposition 3.3]{fukuda2} depends on \cite[Theorem 5.1]{kawamata}. 
Therefore, the proof of Corollary \ref{fuku-main} in this paper is the first rigorous 
proof of Fukuda's important result (cf.~\cite[Proposition 3.3]{fukuda2}). 
Another purpose of this paper 
is to show how to use Shokurov's ideas:~{\em{b-divisors}}, 
{\em{saturation of linear systems}}, various kinds of {\em{adjunction}}, 
and so on, by reproving some known results by our new formulation. 
We recommend this paper as Chapter $8\frac{1}{2}$ of the 
book:~Flips for $3$-folds and $4$-folds. 
We note that this paper is a complement of the 
paper \cite{fujino-kawamata}. 
We do not use the powerful new method developed in 
\cite{ambroQ}, 
\cite{fujino-on}, \cite{fujino-effective}, \cite{fujino-intro}, \cite{fujino-non}, 
\cite{fujino-fundamental}, and \cite{fujino-book}. 
For related topics and applications, see \cite{fujino-finite}, 
\cite[6.~Applications]{gongyo} and 
\cite{cacciola}. 

Let us explain the motivation for our formulation. 

\begin{say}[Motivation]
Let $(X, B)$ be a projective 
klt pair and let $D$ be a nef Cartier divisor on $X$ such that 
$D-(K_X+B)$ is nef and big. 
Then the Kawamata--Shokurov base point free theorem 
means that $|mD|$ is free for every $m\gg 0$. 
Let $f:Y\to X$ be a projective birational morphism 
from a normal projective variety $Y$ such that 
$K_Y+B_Y=f^*(K_X+B)$. 
We note that $f^*D$ is a nef Cartier divisor on $Y$ and that 
$f^*D-(K_Y+B_Y)$ is nef and big. 
It is obvious that $|mf^*D|$ is free for every $m \gg 0$ 
because $|mD|$ is free for every $m \gg 0$. 
In general, we can not directly apply the Kawamata--Shokurov 
base point free theorem 
to $f^*D$ and $(Y, B_Y)$. 
It is because $(Y, B_Y)$ is sub klt but is not always 
klt. 
Note that a $\mathbb Q$-Cartier 
$\mathbb Q$-divisor $L$ on $X$ is nef, big, or semi-ample if and 
only if so is $f^*L$. 
However, the notion of klt is not stable 
under birational pull-backs. 
By adding a {\em{saturation condition}}, 
which is trivially satisfied for klt pairs, we 
can apply the Kawamata--Shokurov base point free theorem for sub klt 
pairs (see Theorem \ref{thm1}). By this new formulation, the base point free theorem 
becomes more flexible and has some important applications. 
\end{say}

\begin{say}[Background]\label{back}  
A key result we need from Ambro's papers is \cite[Theorem 0.2]{ambro1}, 
which is a generalization of \cite[\S 4.~Pull-back of $L_{X/Y}^{ss}$]{fujifuji}. 
It originates from {\em{Kawamata's positivity theorem}} in \cite{kawa2} and 
Shokurov's idea on {\em{adjunction}}. For the details, 
see \cite[\S0.~Introduction]{ambro1}. 
The formulation and calculation we borrow from \cite{ambro2} and \cite{ambro3} 
grew out from Shokurov's {\em{saturation of linear systems}}. 
\end{say}

We summarize the contents of this paper. 
In Section \ref{sec2}, we reformulate the Kawamata--Shokurov 
base point free theorem for sub klt pairs with a 
saturation condition. To state our theorem, we use the 
notion of b-divisors. It is very useful to discuss linear 
systems with some base conditions. 
In Section \ref{sec-new}, we collect basic properties 
of b-divisors and prove some elementary properties. 
In Section \ref{sec4}, 
we discuss a slight generalization of the main theorem of \cite{kawamata}. 
We need this generalization in Section \ref{sec-exc}. 
The main ingredient of our proof is Ambro's 
formulation of Kodaira's canonical 
bundle formula. By this formula and the 
refined Kawamata--Shokurov base point free theorem obtained 
in Section \ref{sec2}, we can quickly prove Kawamata's 
theorem in \cite{kawamata} and its generalization without appealing the notion of 
generalized normal crossing varieties. 
In Section \ref{re-fuku}, we treat the base point free theorem 
of Reid--Fukuda type. In this case, the saturation condition 
behaves very well for inductive arguments. 
It helps us understand the saturation condition of 
linear systems. 
In Section \ref{sec-variants}, we 
prove some variants of 
base point free theorems. They are 
mainly due to Fukuda \cite{fukuda2}. 
We reformulate them by using b-divisors 
and saturation conditions. Then we 
use Ambro's canonical bundle formula to reduce them to 
the easier case instead of 
proving them directly by X-method. 
In Section \ref{sec-exc}, we generalize the Kawamata--Shokurov 
base point free theorem and Kawamata's main theorem in 
\cite{kawamata} for {\em{pseudo-klt pairs}}. 
Theorem \ref{quasi1}, which is new, is the 
main theorem of this section. 
It will be useful for the study of lc centers (cf.~Theorem \ref{final-thm}). 
 
\begin{ack}
The first version of this paper was written 
in Nagoya in 2005.  It was circulated 
as \lq\lq A remark on the base point free theorem\rq\rq 
(arXiv:math/0508554v1). 
The author was partially supported by The Sumitomo Foundation 
and by the Grant-in-Aid for Young Scientists (A) 
$\sharp$17684001 from 
JSPS when he prepared the first version. 
He revised this paper in Kyoto in 2011. 
He was partially supported by The Inamori Foundation and by the 
Grant-in-Aid for Young Scientists (A) $\sharp$20684001 from JSPS. 
\end{ack}

\begin{notation}
Let $B=\sum b_iB_i$ be a $\mathbb Q$-divisor 
on a normal variety $X$ such that $B_i$ is prime for 
every $i$ and that $B_i\ne B_j$ for $i\ne j$. We denote by 
$$
\ulcorner B\urcorner =\sum \ulcorner b_i\urcorner B_i,\ \  
\llcorner B\lrcorner =\sum \llcorner b_i\lrcorner B_i,  
\ \ {\text{and}}\ \ \{B\}=B-\llcorner B\lrcorner 
$$ 
the {\em{round-up}}, the {\em{round-down}}, and 
the {\em{fractional part}} of $B$. 
Note that we do not use $\mathbb R$-divisors in this paper. 
We make one general 
remark here. 
Since the freeness (or semi-ampleness) 
of a Cartier divisor $D$ on a variety $X$ 
depends only on the linear equivalence 
class of $D$, we can 
freely replace $D$ by a linearly equivalent 
divisor to prove the freeness (or semi-ampleness) 
of $D$. 
\end{notation}

We will work over an algebraically closed field $k$ of 
characteristic zero throughout 
this paper. 

\section{Kawamata--Shokurov base point free theorem revisited}\label{sec2}

Kawamata and Shokurov claimed the following theorem 
for klt pairs, that is, they assumed that $B$ is effective. 
In that case, the condition (2) is trivially satisfied. 
We think that our formulation is useful for some applications. 
If the readers are not familiar with 
the notion of b-divisors, we recommend them to see 
Section \ref{sec-new}. 

\begin{thm}[Base point free theorem]\label{thm1}
Let $(X,B)$ be a sub klt pair, let $\pi:X \to S$ be a proper 
surjective morphism onto a variety $S$ and let $D$ be a $\pi$-nef 
Cartier divisor on $X$. Assume the following conditions{\em{:}} 
\begin{itemize}
\item[(1)] $rD-(K_X+B)$ is nef and big over $S$ for some positive integer 
$r$, and 
\item[(2)] $($Saturation condition$)$ 
there exists a positive integer $j_0$ 
such that $\pi_*\mathcal O_X(\ulcorner \mathbf 
A(X,B)\urcorner +j\overline D)\subseteq  
\pi_*\mathcal O_X(jD)$ for every integer $j\geq j_0$. 
\end{itemize}
Then $mD$ is $\pi$-generated for every $m\gg 0$, that is, 
there exists a positive integer $m_0$ such that for 
every $m\geq m_0$ the natural homomorphism $\pi^*\pi_*\mathcal O_X(mD)
\to \mathcal O_X(mD)$ is surjective. 
\end{thm}

\begin{proof}
The usual proof of the base point free theorem, 
that is, X-method, works without any changes if 
we note Lemma \ref{213}. For the details, see, for example, 
\cite[\S 3-1]{kmm}. See also remarks in \ref{four}. 
\end{proof}

The assumptions in Theorem \ref{thm1} are birational in nature. 
This point is indispensable in Section \ref{sec4}. 
We note that we can assume that $X$ is non-singular and 
$\Supp B$ is a simple normal crossing divisor because the conditions 
(1) and (2) are invariant for birational pull-backs. 
So, it is easy to see that Theorem \ref{thm1} 
is equivalent to the following theorem. 

\begin{thm}
Let $X$ be a non-singular variety and let $B$ be a $\mathbb Q$-divisor 
on $X$ such that $\llcorner B\lrcorner \leq 0$ and 
$\Supp B$ is a simple normal crossing 
divisor. 
Let $\pi:X \to S$ be a projective 
morphism onto a variety $S$ and let $D$ be a $\pi$-nef 
Cartier divisor on $X$. Assume the following conditions{\em{:}} 
\begin{itemize}
\item[(1)] $rD-(K_X+B)$ is nef and big over $S$ for some positive integer 
$r$, and 
\item[(2)] $($Saturation condition$)$ 
there exists a positive integer $j_0$ 
such that $\pi_*\mathcal O_X(\ulcorner \mathbf 
-B\urcorner +j D)\simeq 
\pi_*\mathcal O_X(jD)$ for every integer $j\geq j_0$. 
\end{itemize}
Then $mD$ is $\pi$-generated for every $m\gg 0$. 
\end{thm}

The following example says that 
the original Kawamata--Shokurov base point free theorem 
does not necessarily hold for {\em{sub}} klt pairs. 

\begin{ex}\label{e23}
Let $X=E$ be an elliptic curve. 
We take a Cartier divisor $H$ such that 
$\deg H=0$ and $lH\not \sim 0$ for every $l\in \mathbb Z\setminus \{0\}$. 
In particular, $H$ is nef. We put $B=-P$, where 
$P$ is a closed point of $X$. Then $(X, B)$ is sub klt and 
$H-(K_X+B)$ is ample. However, $H$ is not semi-ample. 
In this case, $H^0(X, \mathcal O_X(\ulcorner \mathbf A(X, B)\urcorner 
+j\overline H))\simeq H^0(X, \mathcal O_X(P+jH))\simeq k$ for every $j$. 
However, $H^0(X, \mathcal O_X(jH))=0$ for all $j$. 
Therefore, the saturation condition in Theorem \ref{thm1} 
does not hold. 
\end{ex}

We note that Koll\'ar's effective base point freeness holds under the same 
assumption as in Theorem \ref{thm1}. 

\begin{thm}[Effective freeness]\label{thm23} 
We use the same notation and assumption as in {\em{Theorem \ref{thm1}}}. 
Then there exists a positive integer $l$, which depends only 
on $\dim X$ and $\max\{r, j_0\}$, such that 
$lD$ is $\pi$-generated, that is, 
$\pi^*\pi_*\mathcal O_X(lD)\to \mathcal O_X(lD)$ is 
surjective.  
\end{thm}
\begin{proof}[Sketch of the proof] We need no new ideas. 
So, we just explain how to modify the arguments in \cite[Section 2]{kollar}. 
From now on, we use the notation in \cite{kollar}. 
In \cite{kollar}, $(X, \Delta)$ is assumed to be klt, that is, 
$(X, \Delta)$ is sub klt and $\Delta$ is effective. 
The effectivity of $\Delta$ implies that $H'$ is $f$-exceptional 
in \cite[(2.1.4.3)]{kollar}. We need this to prove 
$H^0(Y, \mathcal O_Y(f^*N+H'))=H^0(X, \mathcal O_X(N))$ in \cite
[(2.1.6)]{kollar}. It is not difficult to see that 
$0\leq H'\leq \ulcorner \mathbf A (X, \Delta)_Y\urcorner$ in our notation. 
Therefore, it is sufficient to assume the saturation condition (the assumption 
(2) in Theorem \ref{thm1}) in 
the proof of Koll\'ar's effective freeness (see Section 2 in 
\cite{kollar}). 
We make one more remark. 
Applying the argument in the first part of 2.4 in \cite{kollar} to 
$\mathcal O_X(j\overline{D}+\ulcorner \mathbf A(X, B)\urcorner)$ on the 
generic fiber of $\pi:X\to S$ with the saturation condition 
(2) in Theorem \ref{thm1}, we obtain a positive integer 
$l_0$ that depends only on $\dim X$ and $\max\{r, j_0\}$ 
such that $\pi_*\mathcal O_X(l_0D)\ne 0$. 
As explained above, the arguments in Section 2 in \cite{kollar} 
work with only minor modifications in our setting. 
We leave the details for the readers' exercise. 
\end{proof}

\section{B-divisors}\label{sec-new}
\begin{say}[Singularities of pairs] 
Let us recall the notion of singularities of pairs. 
We recommend the readers to see \cite{fujino} for 
more advanced topics on singularities of pairs. 

\begin{defn}[Singularities of pairs] 
Let $X$ be a normal variety and let $B$ be a $\mathbb Q$-divisor on $X$ such 
that $K_X+B$ is $\mathbb Q$-Cartier. 
Let $f:Y\to X$ be a resolution of singularities such that 
$\Exc (f)\cup f^{-1}_*B$ has a simple normal crossing support, 
where $\Exc (f)$ is the exceptional locus of $f$. 
We write 
$$
K_Y=f^*(K_X+B)+\sum a_iA_i. 
$$ 
We note that 
$a_i$ is called the {\em{discrepancy}} 
of $A_i$. 
Then the pair $(X,B)$ is {\em{sub klt}} (resp.~{\em{sub lc}}) if 
$a_i>-1$ (resp.~$a_i\geq -1$) for every $i$. 
The pair $(X,B)$ is {\em{klt}} (resp.~{\em{lc}}) if $(X,B)$ is 
sub klt (resp.~sub lc) and $B$ is effective. In some literature, sub klt (resp.~sub 
lc) is sometimes called klt (resp.~lc). 
Let $(X, B)$ be an lc pair. 
If there exists a resolution $f:Y\to X$ 
such that $\Exc (f)$ and $\Exc (f)\cup f^{-1}_*B$ are 
simple normal crossing divisors on $Y$ and 
$$K_Y=f^*(K_X+B)+\sum a_iA_i$$ with 
$a_i>-1$ for all $f$-exceptional $A_i$'s, then 
$(X, B)$ is called {\em{dlt}}.  
\end{defn}

\begin{rem}\label{sub}
Let $(X,B)$ be a klt (resp.~lc) 
pair and $f:Y\to X$ a proper birational morphism of 
normal varieties. We put $K_Y+B_Y=f^*(K_X+B)$. 
Then $(Y, B_Y)$ is not necessarily klt (resp.~lc) 
but sub klt (resp.~sub lc). 
\end{rem}

Let us recall the definition of {\em{log canonical centers}}.   

\begin{defn}[Log canonical center]\label{77} 
Let $(X, B)$ be a sub lc pair. 
A subvariety $W\subset X$ is called a {\em{log canonical center}} 
or an {\em{lc center}} of $(X, B)$ if there 
is a resolution $f:Y\to X$ such that 
$\Exc (f)\cup \Supp f^{-1}_*B$ 
is a simple normal crossing 
divisor on $Y$ and 
a divisor $E$ with discrepancy $-1$ such that 
$f(E)=W$. 
A log canonical center $W\subset X$ of $(X, B)$ is called 
{\em{exceptional}} if there is 
a unique divisor $E_W$ on $Y$ with 
discrepancy $-1$ such that 
$f(E_W)=W$ and 
$f(E)\cap W=\emptyset$ 
for every other divisor $E\ne E_W$ on $Y$ with 
discrepancy $-1$. 
\end{defn}
\end{say}

\begin{say}[b-divisors]
In this paper, we adopt the notion of {\em{b-divisors}}, which 
was introduced by
Shokurov.
For the details of b-divisors, we recommend the readers to see
\cite[1-B]{ambro2} and
\cite[2.3.2]{corti}. 
The readers can find various examples of b-divisors in \cite{isko}. 

\begin{defn}[b-divisor] 
Let $X$ be a normal variety and let $\Div (X)$ be the free abelian 
group generated by Weil divisors on $X$. 
A {\em{b-divisor}} on $X$ is an element: 
$$
\mathbf D \in \mathbf{Div} (X) =\xprojlim_{Y\to X} \Div (Y),  
$$ 
where the projective limit is taken over all proper birational 
morphisms $f:Y\to X$ of normal varieties, under 
the push forward homomorphism $f_*:\Div (Y)\to \Div (X)$. 
A {\em{$\mathbb Q$-b-divisor}} on $X$ is an element of 
$\mathbf {Div}_{\mathbb Q}(X)=\mathbf {Div}(X)\otimes 
_{\mathbb Z}\mathbb Q$. 
\end{defn}

\begin{defn}[Discrepancy $\mathbb Q$-b-divisor]
Let $X$ be a normal variety and let $B$ be a $\mathbb Q$-divisor on
$X$ such that
$K_X+B$ is $\mathbb Q$-Cartier.
Then the {\em{discrepancy $\mathbb Q$-b-divisor}} of the pair $(X,B)$ is
the $\mathbb Q$-b-divisor
$\mathbf A=\mathbf A(X,B)$ with the trace $\mathbf A_Y$ defined by the
formula:
$$
K_Y=f^*(K_X+B)+\mathbf A_Y,
$$
where $f:Y\to X$ is a proper birational morphism of 
normal varieties.
\end{defn}
\begin{defn}[Cartier closure] 
Let $D$ be a $\mathbb Q$-Cartier $\mathbb Q$-divisor on a normal variety
$X$.
Then the $\mathbb Q$-b-divisor $\overline{D}$ denotes the {\em{Cartier
closure}} of $D$, whose trace on $Y$ is $\overline D_Y=f^*D$, 
where $f:Y\to X$ is a proper birational morphism of
normal varieties.
\end{defn}
\end{say}

\begin{defn}
Let $\mathbf D$ be a $\mathbb Q$-b-divisor on $X$. 
The round up $\ulcorner \mathbf D\urcorner \in \mathbf {Div}(X)$ 
is defined componentwise. 
The restriction of $\mathbf D$ to an open subset $U\subset X$ 
is a well-defined $\mathbb Q$-b-divisor on $U$, denoted 
by $\mathbf D |_{U}$. Then $\mathcal O_X(\mathbf D)$ 
is an $\mathcal O_X$-module whose sections on an open subset 
$U\subset X$ are given by 
$$
H^0(U, \mathcal O_X(\mathbf D))=\{a\in k(X)^{\times}; 
(\overline {(a)}+\mathbf D)|_{U}\geq 0\}
\cup \{0\}, 
$$ 
where $k(X)$ is the function field of $X$. 
Note that $\mathcal O_X(\mathbf D)$ is not necessarily 
coherent. 
\end{defn}

\begin{say}[Basic properties] 
We recall the first basic property of discrepancy $\mathbb Q$-b-divisors. 
We will treat a generalization of Lemma \ref{213} 
for sub lc pairs below. 

\begin{lem}\label{213}
Let $(X, B)$ be a sub klt pair and let $D$ be a Cartier divisor on $X$. 
Let $f:Y\to X$ be a proper surjective morphism from a non-singular
variety $Y$.
We write $K_Y=f^*(K_X+B)+\sum a_i A_i$. 
We assume that $\sum  A_i$
is a simple normal crossing divisor.
Then $$\mathcal O_X(\ulcorner \mathbf A(X,B)\urcorner +j\overline{D})=
f_*\mathcal O_Y(\sum \ulcorner a_i\urcorner A_i)\otimes \mathcal O_X(jD)$$ 
for every integer $j$. 

Let $E$ be an effective divisor on $Y$ such that $E\leq \sum \ulcorner
a_i\urcorner A_i$.
Then $$\pi_*f_*\mathcal O_Y(E+f^*jD)\simeq \pi_*\mathcal O_X(jD)$$ if 
$$\pi_*\mathcal O_X(\ulcorner \mathbf A(X, B)\urcorner +j\overline 
D)\subseteq \pi_*\mathcal O_X(jD),$$ where $\pi:X\to S$ is a proper 
surjective morphism onto a variety $S$.
\end{lem}
\begin{proof}
For the first equality, see \cite[Lemmas 2.3.14 and 2.3.15]{corti} 
or their generalizations:~Lemmas \ref{319} and \ref{320} below. 
Since
$E$ is effective,
$\pi_*\mathcal O_X(jD)\subseteq \pi_*f_*\mathcal O_Y(E+f^*jD)\simeq
\pi_*(f_*\mathcal O_Y(E)\otimes \mathcal O_X(jD))$. By 
the assumption and 
$
E\leq \sum
\ulcorner a_i\urcorner A_i$, 
\begin{align*}
\pi_*(f_*\mathcal O_Y(E)\otimes \mathcal
O_X(jD))&\subseteq
\pi_*(f_*\mathcal O_Y(\sum \ulcorner a_i\urcorner A_i)\otimes \mathcal
O_X(jD))\\&=
\pi_*\mathcal O_X(\ulcorner \mathbf A(X,B)\urcorner +j\overline {D})\\ &\subseteq
\pi_*\mathcal O_X(jD). 
\end{align*}
Therefore, we obtain 
$\pi_*f_*\mathcal O_Y(E+f^*jD)\simeq \pi_*\mathcal O_X(jD)$. 
\end{proof}

We will use Lemma \ref{313} 
in Section \ref{sec4}. The 
vanishing theorem in Lemma \ref{313} is nothing but 
the Kawamata--Viehweg--Nadel vanishing theorem. 

\begin{lem}\label{313} 
Let $X$ be a normal variety and let $B$ be a 
$\mathbb Q$-divisor on $X$ such that 
$K_X+B$ is $\mathbb Q$-Cartier. 
Let $f:Y\to X$ be a proper birational 
morphism from a normal variety $Y$. We put 
$K_Y+B_Y=f^*(K_X+B)$. 
Then $$f_*\mathcal O_Y(\ulcorner \mathbf A(Y, B_Y)\urcorner) 
= \mathcal O_X(\ulcorner \mathbf A(X, B)\urcorner)$$ 
and $$R^if_*\mathcal O_Y(\ulcorner 
\mathbf A(X, B) \urcorner )=0$$ for every $i>0$. 
\end{lem}

\begin{proof}
Let $g:Z\to Y$ be a resolution 
such that $\Exc(g)\cup g^{-1}_*B_Y$ has a simple 
normal crossing support. 
We put $K_Z+B_Z=g^*(K_Y+B_Y)$. Then 
$K_Z+B_Z=h^*(K_X+B)$, where 
$h=f\circ g: Z\to X$. By Lemma \ref{213}, 
$$\mathcal O_Y(\ulcorner \mathbf A(Y, B_Y)\urcorner) 
= g_*\mathcal O_Z(\ulcorner -B_Z\urcorner )$$ and 
$$\mathcal O_X(\ulcorner \mathbf A(X, B)\urcorner) 
= h_*\mathcal O_Z(\ulcorner -B_Z\urcorner ). $$
Therefore, 
$f_* \mathcal O_Y(\ulcorner \mathbf A(Y, B_Y)\urcorner)=  
\mathcal O_X(\ulcorner \mathbf A(X, B)\urcorner)$. 
Since, $-B_Z=K_Z-h^*(K_X+B)$, we have 
$$\ulcorner -B_Z\urcorner=K_Z+\{B_Z\}-h^*(K_X+B).$$  
Therefore, $R^ig_*\mathcal O_Z(\ulcorner -B_Z\urcorner)=0$ and 
$R^ih_*\mathcal O_Z(\ulcorner -B_Z\urcorner)=0$ for 
every $i>0$ by the Kawamata--Viehweg vanishing 
theorem. Thus, $$R^if_*\mathcal O_Y(\ulcorner \mathbf A(Y, 
B_Y)\urcorner)=0$$ for every $i>0$ by 
Leray's spectral sequence. 
\end{proof}

\begin{rem}\label{214} 
We use the same notation as in Remark \ref{sub}. 
Let $(X, B)$ be a klt pair. 
Let $D$ be a Cartier divisor on $X$ and 
let $\pi:X\to S$ be a proper morphism onto a variety $S$. 
We put $p=\pi\circ f:Y\to S$. 
Then $p_*\mathcal O_Y(jf^*D)\simeq \pi_*\mathcal O_X(jD)\simeq 
p_*\mathcal O_Y(\ulcorner \mathbf A(Y, B_Y)\urcorner +j{\overline{f^*D}})$ for 
every integer $j$. 
It is because $f_*\mathcal O_Y(\ulcorner 
\mathbf A(Y, B_Y)\urcorner )= 
\mathcal O_X(\ulcorner \mathbf A(X, B)\urcorner )\simeq \mathcal O_X$ 
by Lemma \ref{313}. 
\end{rem}

We make a brief comment on the {\em{multiplier ideal sheaf}}. 
\begin{rem}[Multiplier ideal sheaf] 
Let $D$ be an effective $\mathbb Q$-divisor on 
a non-singular variety $X$. Then $\mathcal O_X(\ulcorner \mathbf A 
(X, D)\urcorner )$ is nothing but the {\em{multiplier 
ideal sheaf}} $\mathcal J (X, D)\subseteq \mathcal O_X$ of $D$ 
on $X$. See \cite[Definition 9.2.1]{lazarsfeld}. 
More generally, let $X$ be a normal variety and let $\Delta$ be a $\mathbb Q$-divisor on 
$X$ such that $K_X+\Delta$ is $\mathbb Q$-Cartier. 
Let $D$ be a $\mathbb Q$-Cartier $\mathbb Q$-divisor on $X$. 
Then $\mathcal O_X(\ulcorner \mathbf A(X, \Delta +D)\urcorner )
=\mathcal J((X,\Delta); D)$, where the right hand side is the 
{\em{multiplier ideal sheaf}} defined (but not investigated) 
in \cite[Definition 9.3.56]{lazarsfeld}. 
In general, $\mathcal O_X(\ulcorner \mathbf A(X, \Delta+D)\urcorner )$ is a fractional 
ideal of $k(X)$. 
\end{rem}
\end{say} 

\begin{say}[Remarks on Theorem \ref{thm1}]\label{four} 
The following four remarks help us understand 
Theorem \ref{thm1}. 
\begin{rem}[Non-vanishing theorem]
By Shokurov's 
non-vanishing theorem (see \cite[Theorem 2-1-1]{kmm}), 
we have 
that $\pi_*\mathcal O_X(\ulcorner \mathbf 
A(X,B)\urcorner +j\overline D)\ne 0$ for every $j\gg 0$. Thus we 
have $\pi_*\mathcal O_X(jD)\ne 0$ for every $j\gg 0$ by 
the condition (2) in Theorem \ref{thm1}. 
\end{rem}
\begin{rem}
We know that $\ulcorner \mathbf 
A(X,B)\urcorner\geq 0$ since $(X,B)$ is sub klt. 
Therefore, $\pi_*\mathcal O_X(jD)\subseteq 
\pi_*\mathcal O_X(\ulcorner \mathbf 
A(X,B)\urcorner +j\overline D)$. This implies that 
$\pi_*\mathcal O_X(jD)\simeq 
\pi_*\mathcal O_X(\ulcorner \mathbf 
A(X,B)\urcorner +j\overline D)$ for $j\geq j_0$ by the condition (2) 
in Theorem \ref{thm1}. 
\end{rem}
\begin{rem}\label{3131}
If the pair $(X,B)$ is klt, then 
$\ulcorner \mathbf 
A(X,B)\urcorner$ is effective and exceptional over $X$. 
In this case, it is obvious that $\pi_*\mathcal O_X(jD)=
\pi_*\mathcal O_X(\ulcorner \mathbf 
A(X,B)\urcorner +j\overline D)$. 
\end{rem}
\begin{rem}
The condition (2) in Theorem \ref{thm1} 
is a very elementary case of {\em{saturation of linear 
systems}}. See \cite[2.3.3]{corti} and \cite[1-D]{ambro2}. 
\end{rem}
\end{say}

\begin{say} We introduce the notion 
of {\em{non-klt $\mathbb Q$-b-divisor}}, which 
is trivial for sub klt pairs. We will use this in Section \ref{re-fuku}. 

\begin{defn}[Non-klt $\mathbb Q$-b-divisor] 
Let $X$ be a normal variety and let $B$ be a $\mathbb Q$-divisor 
on $X$ such that $K_X+B$ is $\mathbb Q$-Cartier. 
Then the {\em{non-klt $\mathbb Q$-b-divisor}} of the 
pair $(X, B)$ is the $\mathbb Q$-b-divisor $\mathbf N=\mathbf N(X, B)$ 
with the trace $\mathbf N_Y=\sum_{a_i\leq -1}a_iA_i$ for 
$$
K_Y=f^*(K_X+B)+\sum a_iA_i, 
$$ 
where $f:Y\to X$ is a proper birational morphism of normal 
varieties. It is easy to see that $\mathbf N(X, B)$ is a well-defined 
$\mathbb Q$-b-divisor. 
We put $\mathbf A^{\ast} (X, B)=\mathbf A(X, B)-\mathbf N(X, 
B)$. Of course, $\mathbf A^{\ast}(X, B)$ is a well-defined 
$\mathbb Q$-b-divisor and $\ulcorner 
\mathbf A^{\ast}(X, B)\urcorner\geq 0$. 
If $(X, B)$ is sub klt, then $\mathbf N(X, B)=0$ and 
$\mathbf A(X, B)=\mathbf A^{\ast} (X, B)$. 
\end{defn}

The next lemma is a generalization of Lemma \ref{213}. 

\begin{lem}\label{319}
Let $(X, B)$ be a sub lc pair and 
let $f:Y\to X$ be a resolution such that 
$\Exc (f)\cup \Supp f^{-1}_*B$ is a simple 
normal crossing divisor on $Y$. 
We write $K_Y=f^*(K_X+B)+\sum a_iA_i$. Then 
$$\mathcal O_X(\ulcorner \mathbf A^{\ast}(X, B)\urcorner)=f_*\mathcal O
_Y(\sum _{a_i\ne -1}\ulcorner a_i\urcorner A_i). $$In particular, 
$\mathcal O_X(\ulcorner \mathbf A^*(X, B)\urcorner)$ is 
a coherent $\mathcal O_X$-module. If $(X, B)$ is lc, then 
$\mathcal O_X(\ulcorner \mathbf A^{\ast}(X, B)\urcorner)\simeq 
\mathcal O_X$. 

Let $D$ be a Cartier divisor on $X$ and let $E$ be an effective divisor on $Y$ such 
that $E\leq \sum _{a_i\ne -1} \ulcorner a_i\urcorner A_i$. 
Then $$\pi_*f_*\mathcal O_Y(E+f^*jD)\simeq \pi_*\mathcal 
O_X(D)$$ if $$\pi_*\mathcal O_X(\ulcorner \mathbf A^{\ast}(X, B)\urcorner 
+j\overline D)\subseteq \pi_*\mathcal O_X(jD), $$
where $\pi:X\to S$ is a proper morphism onto a variety $S$.   
\end{lem}

\begin{proof}
By definition, $\mathbf A^{\ast}(X, B)_Y=\sum _{a_i\ne -1}a_iA_i$. 
If $g:Y'\to Y$ is a proper birational 
morphism 
from a normal variety $Y'$, then 
$$
\ulcorner \mathbf A^{\ast}(X, B)_{Y'}\urcorner=g^*
\ulcorner \mathbf A^{\ast}(X,B)_Y\urcorner+F, 
$$ 
where $F$ is a $g$-exceptional effective divisor, by Lemma 
\ref{320} below. 
This implies $f_*\mathcal O_Y
(\ulcorner \mathbf A^{\ast}(X, B)_{Y}\urcorner)=
f'_*\mathcal O_{Y'}(\ulcorner \mathbf A^{\ast}(X, B)_{Y'}\urcorner)$, 
where $f'=f\circ g$, from which it follows 
that 
$\mathcal O_X(\ulcorner \mathbf A^{\ast}(X, B)\urcorner)=f_*\mathcal O
_Y(\sum _{a_i\ne -1}\ulcorner a_i\urcorner A_i)$ is a 
coherent $\mathcal O_X$-module. 
The latter statement is easy to check. 
\end{proof}

\begin{lem}\label{320} 
Let $(X, B)$ be a sub lc pair and $f:Y\to X$ be a resolution 
as in {\em{Lemma \ref{319}}}. We 
consider the 
$\mathbb Q$-b-divisor 
$\mathbf A^{\ast}=\mathbf A^{\ast}(X, B)=
\mathbf A(X, B)-\mathbf N(X, B)$. 
If $Y'$ is a normal variety and $g:Y'\to Y$ is a proper 
birational morphism, then 
$$
\ulcorner \mathbf A^{\ast}_{Y'}\urcorner=g^*
\ulcorner \mathbf A^{\ast}_Y\urcorner+F, 
$$
where $F$ is a $g$-exceptional effective divisor. 
\end{lem}
\begin{proof}
By definition, we have $K_Y=f^*(K_X+B)+\mathbf A_Y$. 
Therefore, we may write, 
\begin{align*}
K_{Y'}&=g^*f^*(K_X+B)+\mathbf A_{Y'}\\ 
&=g^*(K_Y-\mathbf A_Y)+\mathbf A_{Y'}\\ 
&=g^*(K_Y+\{-\mathbf A^{\ast}_Y\}-\mathbf N_Y+\llcorner 
-\mathbf A^{\ast}_Y\lrcorner)+\mathbf A_{Y'}\\ 
&=g^*(K_Y+\{-\mathbf A^{\ast}_Y\}-\mathbf N_Y)+\mathbf A_{Y'}
-g^*\ulcorner \mathbf A^{\ast}_Y\urcorner. 
\end{align*}
We note that $(Y, \{-\mathbf A^{\ast}_Y\}-\mathbf N_Y)$ is lc and 
that the set of lc centers of $(Y, \{-\mathbf A^{\ast}_Y\}-\mathbf N_Y)$ 
coincides with 
that of $(Y, -\mathbf A^{\ast}_Y-\mathbf N_Y)=(Y, -\mathbf A_Y)$. 
Therefore, the round-up of $\mathbf A_{Y'}-g^*\ulcorner 
\mathbf A^{\ast}_Y\urcorner -\mathbf N_{Y'}$ 
is effective and $g$-exceptional. Thus, we can write 
$
\ulcorner \mathbf A^{\ast}_{Y'}\urcorner=g^*
\ulcorner \mathbf A^{\ast}_Y\urcorner+F, 
$
where $F$ is a $g$-exceptional effective divisor. 
\end{proof}
The next lemma is obvious by Lemma \ref{319}. 

\begin{lem}\label{lem321}
Let $(X, B)$ be a sub lc pair and 
let $f:Y\to X$ be a proper birational morphism from a 
normal variety $Y$. We put 
$K_Y+B_Y=f^*(K_X+B)$. Then 
$f_*\mathcal O_Y(\ulcorner \mathbf A^{\ast}(Y, B_Y)\urcorner)
=\mathcal O_X(\ulcorner \mathbf A^{\ast}(X, B)\urcorner)$. 
\end{lem}
\end{say}

\section{Base point free theorem; nef and abundant case}\label{sec4}
We recall the definition of {\em{abundant}} divisors, which are called 
{\em{good}} divisors in \cite{kawamata}. See \cite[\S 6-1]{kmm}. 

\begin{defn}[Abundant divisor] 
Let $X$ be a complete normal 
variety and let $D$ be a $\mathbb Q$-Cartier nef $\mathbb Q$-divisor on $X$. 
We define the {\em{numerical Iitaka dimension}} to be 
$$
\nu(X,D)=\max \{e; D^e\not\equiv 0\}. 
$$ 
This means that $D^{e'}\cdot S=0$ for any 
$e'$-dimensional subvarieties $S$ of $X$ with 
$e'>e$ and there exists an $e$-dimensional 
subvariety $T$ of $X$ such that 
$D^e \cdot T>0$. 
Then it is easy to see that $\kappa (X,D)\leq \nu(X, D)$, where 
$\kappa (X,D)$ denotes {\em{Iitaka's $D$-dimension}}. 
A nef $\mathbb Q$-divisor $D$ is said to be {\em{abundant}} if the equality 
$\kappa (X, D)=\nu(X, D)$ holds. 
Let $\pi:X\to S$ be a proper surjective morphism of normal varieties and 
let $D$ be a $\mathbb Q$-Cartier $\mathbb Q$-divisor on $X$. 
Then $D$ is said to be $\pi$-abundant if $D|_{X_\eta}$ is abundant, where 
$X_\eta$ is the generic fiber of $\pi$.   
\end{defn}

The next theorem is 
the main theorem of \cite{kawamata}. 
For the relative statement, see 
\cite[Theorem 5]{nakayama}. 
We reduced 
Theorem \ref{thm2} to Theorem \ref{thm1} 
by using Ambro's results in \cite{ambro1} and \cite{ambro3}, 
which is the main theme of \cite{fujino-kawamata}. 
For the details, see \cite[Section 2]{fujino-kawamata}. 

\begin{thm}[{cf.~\cite[Theorem 6-1-11]{kmm}}]\label{thm2}
Let $(X,B)$ be a klt pair, let $\pi:X \to S$ be a proper 
morphism onto a variety $S$. Assume the following 
conditions{\em{:}} 
\begin{itemize}
\item[(a)] $H$ is a $\pi$-nef $\mathbb Q$-Cartier $\mathbb Q$-divisor on $X$, 
\item[(b)] $H-(K_X+B)$ is $\pi$-nef and $\pi$-abundant, and 
\item[(c)] $\kappa (X_\eta, (aH-(K_X+B))_\eta)\geq 0$ 
and $\nu(X_\eta, (aH-(K_X+B))_\eta)=\nu(X_\eta, (H-(K_X+B))_\eta)$ for 
some $a\in \mathbb Q$ with $a>1$, 
where $\eta$ is the generic point of $S$. 
\end{itemize} 
Then $H$ is $\pi$-semi-ample. 
\end{thm}

We recall the definition of the Iitaka fibrations in this paper 
before we state the main theorem of this section. 
\begin{defn}[Iitaka fibration]
Let $\pi:X\to S$ be a proper surjective morphism of normal varieties.
Let $D$ be a $\mathbb Q$-Cartier $\mathbb Q$-Weil divisor on $X$ such that
$\kappa(X_{\eta}, D_{\eta})\geq 0$, where $\eta$ is the generic point of $S$.
Let $X\dashrightarrow W$ be the rational map over $S$ induced by $\pi^*\pi_*\mathcal
O_X(mD)\to \mathcal O_X(mD)$ for a sufficiently large and divisible
integer $m$.
We consider a projective surjective morphism $f:Y\to Z$ of non-singular
varieties that is birational to $X\dashrightarrow W$. We call
$f:Y\to Z$ the {\em{Iitaka fibration}} with respect to $D$ over $S$.
\end{defn}

Theorem \ref{thm4} is a slight generalization of 
Theorem \ref{thm2}. 
It will be used in the proof of Theorem \ref{quasi1}. 

\begin{thm}\label{thm4}
Let $(X,B)$ be a sub klt pair, let $\pi:X \to S$ be a proper 
morphism onto a variety $S$. Assume the following 
conditions{\em{:}} 
\begin{itemize}
\item[(a)] $H$ is a $\pi$-nef $\mathbb Q$-Cartier $\mathbb Q$-divisor on $X$, 
\item[(b)] $H-(K_X+B)$ is $\pi$-nef and $\pi$-abundant, 
\item[(c)] $\kappa (X_\eta, (aH-(K_X+B))_\eta)\geq 0$ 
and $\nu(X_\eta, (aH-(K_X+B))_\eta)=\nu(X_\eta, (H-(K_X+B))_\eta)$ for 
some $a\in \mathbb Q$ with $a>1$, 
where $\eta$ is the generic point of $S$,  
\item[(d)] let $f:Y\to Z$ be the Iitaka fibration with respect to $H-(K_X+B)$ over
$S$. We assume that there exists a proper birational morphism $\mu: Y\to X$
and put $K_Y+B_Y=\mu ^*(K_X+B)$. In this setting, we assume 
$\rank f_*\mathcal O_Y(\ulcorner \mathbf A(Y, B_Y)\urcorner)=1$, and 
\item[(e)] $($Saturation condition$)$
there exist positive integers $b$ and $j_0$ such that
$bH$ is Cartier and $\pi_*\mathcal O_X(\ulcorner
\mathbf A(X, B)\urcorner+jb\overline H)\subseteq \pi_*
\mathcal O_X(jbH)$ for every positive integer $j\geq j_0$. 
\end{itemize}
Then $H$ is $\pi$-semi-ample. 
\end{thm}
\begin{proof}
The proof of Theorem \ref{thm2} in 
\cite[Section 2]{fujino-kawamata} 
works without any changes. 
We note that the condition (d) implies \cite[Lemma 2.3]{fujino-kawamata} 
and that we can use the condition (e) in the proof of 
\cite[Lemma 2.4]{fujino-kawamata}. 
\end{proof}
\begin{rem}
We note that $\rank f_*\mathcal O_Y(\ulcorner \mathbf A(Y, B_Y)
\urcorner )$ is a birational invariant for $f:Y\to Z$ by Lemma 
\ref{313}. 
\end{rem}

\begin{rem}\label{4646}
If $(X,B)$ is klt and $bH$ is Cartier, then it is obvious that 
$\pi_*\mathcal O_X(\ulcorner \mathbf A(X, B)\urcorner +jb\overline 
H)\simeq \pi_* \mathcal O_X(jbH)$ for every positive 
integer $j$ (see Remark \ref{3131}). 
\end{rem}

\begin{rem}We can easily generalize Theorem \ref{thm4} 
to varieties in class $\mathcal C$ by suitable modifications. 
For details, see \cite[Section 4]{fujino-kawamata}. 
\end{rem}

The following examples help us to understand 
the condition (d). 

\begin{ex}
Let $X=E$ be an elliptic curve and $P\in X$ a closed 
point. Take a general member $P_1+P_2+P_3\in 
|3P|$. We put $B=\frac{1}{3}(P_1+P_2+P_3)-P$. 
Then $(X, B)$ is sub klt and $K_X+B\sim _{\mathbb Q}0$. 
In this case, $\mathcal O_X(\ulcorner \mathbf A(X, B)\urcorner)
\simeq \mathcal O_X(P)$ and $H^0(X, \mathcal O_X
(\ulcorner \mathbf A(X, B)\urcorner))\simeq k$. 
\end{ex}

\begin{ex}\label{e48}
Let $f:X=\mathbb P_{\mathbb P^1}(\mathcal O_{\mathbb P^1} 
\oplus \mathcal O_{\mathbb P^1}(1))\to Z=\mathbb P^1$ 
be the Hirzebruch surface and let $C$ (resp.~$E$) be the 
positive (resp.~negative) section of $f$. 
We take a general member $B_0\in |5C|$. Note that 
$|5C|$ is a free linear system on $X$. 
We put $B=-\frac{1}{2} E+\frac{1}{2} B_0$ and consider 
the pair $(X, B)$. Then $(X, B)$ is sub klt. We put 
$H=0$. Then $H$ is a nef Cartier divisor on $X$ and 
$aH-(K_X+B)\sim_{\mathbb Q}\frac{1}{2} F$ for every rational 
number $a$, where $F$ is a fiber of $f$. 
Therefor, $aH-(K_X+B)$ is nef and abundant for every rational 
number $a$. In this case, $\mathcal O_X(\ulcorner 
\mathbf A(X, B)\urcorner)\simeq \mathcal O_X(E)$. So, 
we have 
\begin{align*}
H^0(X, \mathcal O_X(\ulcorner 
\mathbf A(X, B)\urcorner+j\overline H))&\simeq H^0(X, 
\mathcal O_X(E))\simeq k
\\ &\simeq H^0(X, \mathcal O_X)
\simeq H^0(X, \mathcal O_X(jH))
\end{align*} for every integer $j$. 
Therefore, $\pi:X\to \Spec \,k$, $H$, and $(X, B)$ satisfy the 
conditions (a), (b), (c), and (e) in Theorem \ref{thm4}. 
However, (d) is not satisfied. 
In our case, it is easy to see that $f:X\to Z$ is the Iitaka 
fibration with 
respect to $H-(K_X+B)$. 
Since $f_*\mathcal O_X(\ulcorner \mathbf A(X, B)\urcorner)
\simeq f_*\mathcal O_X(E)$, we have 
$\rank f_*\mathcal O_X(\ulcorner 
\mathbf A(X, B)\urcorner)=2$.   
\end{ex}

\begin{rem}\label{rem49}
In Theorem \ref{thm4}, the assumptions (a), (b), 
(c) are the same as in Theorem \ref{thm2}. 
The condition (e) is indispensable 
by Example \ref{e23} for sub klt pairs. 
By using the non-vanishing theorem for generalized 
normal crossing varieties in \cite[Theorem 5.1]{kawamata}, 
which is the hardest part to prove in \cite{kawamata}, 
the semi-ampleness of $H$ seems to follow from the 
conditions (a), (b), (c), and (e). 
However, we need (d) to apply Ambro's canonical 
bundle formula to the Iitaka fibration $f:Y\to Z$. 
See, for example, \cite[Section 3]{fujino-kawamata}. 
Unfortunately, 
as we saw in Example \ref{e48}, the condition (d) does not follow from 
the other assumptions. Anyway, the condition (d) is automatically 
satisfied if $(X, B)$ is klt (see \cite[Lemma 2.3]{fujino-kawamata}). 
\end{rem}

\begin{say}[Examples] 
The following two examples show that 
the effective version of Theorem \ref{thm2} does not 
necessarily hold. 
The first one is an obvious example. 

\begin{ex}
Let $X=E$ be an elliptic curve and let $m$ be an arbitrary 
positive integer. Then there is a Cartier divisor 
$H$ on $X$ such that $mH\sim 0$ and $lH\not\sim 0$ for $0<l<m$. 
Therefore, the effective version of Theorem \ref{thm2} does 
not necessarily hold. 
\end{ex}

The next one shows the reason why 
Theorem \ref{thm23} does not imply the 
effective version of Theorem \ref{thm2}.  

\begin{ex}
Let $E$ be an elliptic curve and $G=\mathbb Z/m\mathbb Z=\langle 
\zeta\rangle$, 
where $\zeta$ is a primitive 
$m$-th root of unity. We take an $m$-torsion point $a\in E$. 
The cyclic group $G$ acts on $E\times 
\mathbb P^1$ as follows:
$$E\times \mathbb P^1\ni 
(x, [X_0:X_1])\mapsto (x+a, [\zeta X_0:X_1])\in 
E\times \mathbb P^1. $$ 
We put $X=(E\times \mathbb P^1)/G$. 
Then $X$ has a structure of elliptic surface $p:X\to \mathbb P^1$. 
In this setting, $K_X=p^*(K_{\mathbb P^1}+\frac{m-1}{m} [0] 
+\frac{m-1}{m}[\infty])$. We put $H=p^{-1}(0)_{\red}$. Then 
$H$ is a Cartier divisor on $X$. It is easy to see that 
$H$ is nef and $H-K_X$ is nef and abundant. 
Moreover, $\kappa (X, aH-K_X)=\nu (X, aH-K_X)=1$ for 
every rational number $a>0$. It is obvious that 
$|mH|$ is free. However, $|lH|$ is not free for $0<l<m$. 
Thus, the effective version of Theorem \ref{thm2} 
does not hold. 
\end{ex}
\end{say} 

\section{Base point free theorem of Reid--Fukuda type}\label{re-fuku}

The following theorem is a reformulation of the main theorem 
of \cite{fuji-tokyo}. 

\begin{thm}[Base point free theorem of Reid--Fukuda type]\label{rftype} 
Let $X$ be a non-singular variety and let $B$ be a $\mathbb Q$-divisor 
on $X$ such that $\Supp B$ is a simple normal crossing divisor 
and $(X, B)$ is sub lc. 
Let $\pi:X\to S$ be a proper morphism onto a 
variety $S$ and let $D$ be a $\pi$-nef Cartier divisor on $X$. 
Assume the following conditions{\em{:}} 
\begin{itemize}
\item[(1)] $rD-(K_X+B)$ is nef and log big over $S$ for some 
positive integer $r$, and 
\item[(2)] $($Saturation condition$)$ 
there exists a positive integer $j_0$ such that 
$\pi_*\mathcal O_X(\ulcorner \mathbf A^{\ast}(X, B)\urcorner
+j\overline D)\subseteq 
\pi_*\mathcal O_X(jD)$ for 
every integer $j\geq j_0$. 
\end{itemize}
Then $mD$ is $\pi$-generated for every $m\gg 0$, that is, there 
exists a positive integer $m_0$ such that 
for every $m\geq m_0$ the natural 
homomorphism $\pi^*\pi_*\mathcal O_X(mD)\to \mathcal O_X(mD)$ 
is surjective. 
\end{thm}

Let us recall the definition of {\em{nef and log big divisors}} 
on sub lc pairs. 
\begin{defn}
Let $(X, B)$ be a sub lc pair and let $\pi:X\to S$ be a proper 
morphism onto a variety $S$. 
Let $\mathcal L$ be a line bundle on $X$. We say that 
$\mathcal L$ is {\em{nef and log big}} over $S$ if and only if 
$\mathcal L$ is $\pi$-nef and $\pi$-big and the 
restriction $\mathcal L|_{W}$ is big over 
$\pi(W)$ for every lc center $W$ of the pair $(X, B)$. 
A $\mathbb Q$-Cartier $\mathbb Q$-divisor 
$H$ on $X$ is said to be nef and log big over $S$ if and 
only if so is $\mathcal O_X(cH)$, where 
$c$ is a positive integer such that $cH$ is Cartier. 
\end{defn}

\begin{proof}[Proof of {\em{Theorem \ref{rftype}}}] 
We write $B=T+B_+-B_-$ such that $T$, $B_+$, and 
$B_-$ are effective divisors, they have no common irreducible 
components, $\llcorner B_+\lrcorner =0$, and $\llcorner T\lrcorner =T$. 
If $T=0$, then $(X, B)$ is sub klt. So, theorem 
follows from Theorem \ref{thm1}. 
Thus, we assume $T\ne 0$. 
Let $T_0$ be an irreducible component of $T$. 
If $m\geq r$, then 
\begin{align*}
mD+\ulcorner B_-\urcorner -T_0-(K_X+B+\ulcorner 
B_-\urcorner-T_0)=mD-(K_X+B)
\end{align*} is nef and log big over $S$ for 
the pair $(X, B+\ulcorner B_-\urcorner -T_0)$. 
We note that $B+\ulcorner B_-\urcorner -T_0$ is effective. 
Therefore, $R^1\pi_*\mathcal O_X(\ulcorner B_-\urcorner 
-T_0+mD)=0$ for 
$m\geq r$ by the vanishing theorem:~Lemma \ref{5vani}. 
Thus, we obtain the following 
commutative 
diagram for $m\geq \max \{r, j_0\}$: 

$$
\begin{CD}
\pi_*\mathcal O_X(\ulcorner B_-\urcorner +mD)
@>>> \pi_*\mathcal O_{T_0}(\ulcorner 
B_-|_{T_0}\urcorner +mD|_{T_0})@>>>0\\
@AA\text{$\cong$}A @AA\text{$\iota$}A\\
\pi_*\mathcal O_X(mD)
@>>>\pi_*\mathcal O_{T_0}(mD|_{T_0}) . 
\end{CD}
$$
Here, we used 
\begin{align*}
\pi_*\mathcal O_X(mD)&\subseteq \pi_*\mathcal O_X(\ulcorner B_-
\urcorner +mD)\\ &\simeq \pi_*\mathcal O_X(\ulcorner 
\mathbf A^{\ast}(X, B)\urcorner +m\overline D)\\ &\subseteq 
\pi_*\mathcal O_X(mD)
\end{align*} for 
$m\geq j_0$ (see Lemma \ref{319}). 
We put $K_{T_0}+B_{T_0}=(K_X+B)|_{T_0}$ and $D_{T_0}=D|_{T_0}$. 
Then $(T_0, B_{T_0})$ is sub lc and it is easy to 
see that $rD_{T_0}-(K_{T_0}+B_{T_0})$ is nef and log big 
over $\pi(T_0)$. 
It is obvious that $T_0$ is non-singular and 
$\Supp B_{T_0}$ is a simple normal crossing divisor. 
We note that $\pi_*\mathcal O_{T_0}(\ulcorner \mathbf A^{\ast}
(T_0, B_{T_0})\urcorner +j\overline {D_{T_0}})\simeq  
\pi_*\mathcal O_{T_0}(jD_{T_0})$ for 
every $j\geq \max \{r, j_0\}$ follows from the above diagram, 
that is, the natural inclusion $\iota$ is isomorphism 
for $m\geq \max \{r, j_0\}$. 
By induction, $mD_{T_0}$ is $\pi$-generated for every 
$m\gg 0$. We can apply the same argument to every irreducible 
component of $T$. 
Therefore, the relative base locus of $mD$ is disjoint 
from $T$ for every $m\gg 0$ since the restriction 
map $\pi_*\mathcal O_X(mD)\to 
\pi_*\mathcal O_{T_0}(mD_{T_0})$ is surjective 
for every irreducible component $T_0$ of $T$. 
By the same arguments as in \cite[Proof of 
Theorem 3]{fukuda}, which is a variant 
of X-method, work without any changes. 
So, we obtain that $mD$ is $\pi$-generated for 
every $m\gg 0$. 
\end{proof}

The following vanishing theorem 
was already used in the proof of Theorem \ref{rftype}. 
The proof is an easy exercise by the induction on $\dim X$ and 
the number of the irreducible components of $\llcorner \Delta\lrcorner$. 

\begin{lem}\label{5vani} 
Let $\pi:X\to S$ be a proper morphism from a non-singular variety $X$. 
Let $\Delta =\sum d_i\Delta_i$ be a sum of 
distinct prime divisors such that 
$\Supp \Delta$ is a simple normal crossing 
divisor and $d_i$ is a rational number with $0\leq d_i\leq 1$ 
for every $i$. Let $D$ be a Cartier divisor on $X$. 
Assume that $D-(K_X+\Delta)$ is nef and log big over 
$S$ for the pair $(X, \Delta)$. Then 
$R^i\pi_*\mathcal O_X(D)=0$ for 
every $i>0$. 
\end{lem}

As in Theorem \ref{thm23}, effective freeness holds under the same 
assumption as in Theorem \ref{rftype}. 
 
\begin{thm}[Effective freeness]\label{52}  
We use the same notation and assumption as in 
{\em{Theorem \ref{rftype}}}. 
Then there exists a positive integer $l$, which depends only 
on $\dim X$ and $\max\{r, j_0\}$, such that 
$lD$ is $\pi$-generated, that is, 
$\pi^*\pi_*\mathcal O_X(lD)\to \mathcal O_X(lD)$ is 
surjective.  
\end{thm}
\begin{proof}[Sketch of the proof] 
If $(X, B)$ is sub klt, then this theorem is nothing but Theorem \ref{thm23}. 
So, we can assume that $(X, B)$ is not sub klt. 
In this case, the arguments in \cite[\S4]{fukuda} 
work with only minor modifications. 
From now on, we use the notation in \cite[\S4]{fukuda}. 
By minor modifications, the proof in \cite[\S4]{fukuda} 
works under the following weaker assumptions:~$X$ is non-singular 
and $\Delta$ is a $\mathbb Q$-divisor on $X$ such that 
$\Supp \Delta$ is a simple normal crossing divisor and 
$(X, \Delta)$ is sub lc.  
In \cite[Claim 5]{fukuda}, $E_i$ is $f$-exceptional. 
In our setting, it is not true. However, 
$0\leq \sum _{cb_i-e_i+p_i<0}\ulcorner 
-(cb_i-e_i+p_i)\urcorner E_i\leq \ulcorner 
\mathbf A^{\ast}(X, \Delta)_Y\urcorner$, 
which always holds even when $\Delta$ is not effective, 
is sufficient for us. 
It is because we can use the saturation condition (2) in Theorem 
\ref{rftype}. 
We leave the details for the readers' exercise since 
all we have to do is to repeat the arguments in \cite[Section 2]{kollar} and 
\cite[\S4]{fukuda}. 
\end{proof}

The final statement in this section is the (effective) 
base point free theorem 
of Reid--Fukuda type for dlt pairs. 

\begin{cor}
Let $(X, B)$ be a dlt pair and let $\pi:X\to S$ be a proper 
morphism onto a variety $S$. 
Let $D$ be a $\pi$-nef Cartier divisor 
on $X$. 
Assume that $rD-(K_X+B)$ is nef and log big over 
$S$ for some positive integer $r$. 
Then there exists a positive integer 
$m_0$ such that $mD$ is $\pi$-generated for 
every $m\geq m_0$ and we can find a positive integer 
$l$, which depends only on $\dim X$ and $r$, such that 
$lD$ is $\pi$-generated. 
\end{cor}
\begin{proof}
Let $f:Y\to X$ be a resolution such that 
$\Exc(f)$ and $\Exc(f)\cup \Supp f^{-1}_*B$ are 
simple normal crossing divisors, 
$K_Y+B_Y=f^*(K_X+B)$, and $f$ is an isomorphism over all 
the generic points of lc 
centers of the pair $(X, B)$. 
Then $(Y, B_Y)$ is sub lc, and $rD_Y-(K_Y+B_Y)$ is nef and 
log big over $S$, where $D_Y=f^*D$. 
Since $\ulcorner \mathbf A^{\ast}(X, B)\urcorner$ is effective and 
exceptional over $X$, 
$p_*\mathcal O_Y(\ulcorner \mathbf A^{\ast}(Y, B_Y)\urcorner 
+j\overline {D_Y})\simeq p_*\mathcal O_Y(jD_Y)$ for 
every $j$, where $p=\pi\circ f$. So, we can apply Theorems 
\ref{rftype} and \ref{52} to $D_Y$ and $(Y, B_Y)$. 
We finish the proof. 
\end{proof}

For the (effective) base point freeness for lc pairs, 
see \cite{fujino-effective}, 
\cite[Theorem 1.2]{fujino-non}, \cite[Theorem 13.1]{fujino-fundamental}, 
and \cite[3.3.1 Base Point Free Theorem]{fujino-book}. 

\section{Variants of base point free theorems due to 
Fukuda}\label{sec-variants}

The starting point of this section is a slight generalization 
of Theorem \ref{thm1}. 
It is essentially the same as \cite[Theorem 3]{fukuda}. 

\begin{thm}\label{bigfukuda}
Let $X$ be a non-singular variety and let $B$ be a $\mathbb Q$-divisor 
on $X$ such that $(X, B)$ is sub lc and 
$\Supp B$ is a simple normal crossing 
divisor. 
Let $\pi:X \to S$ be a proper 
morphism onto a variety $S$ and let $H$ be a $\pi$-nef 
$\mathbb Q$-Cartier $\mathbb Q$-divisor 
on $X$. Assume the following conditions{\em{:}} 
\begin{itemize}
\item[(1)] $H-(K_X+B)$ is nef and big over $S$, and 
\item[(2)] $($Saturation condition$)$ 
there exist positive integers $b$ and $j_0$ 
such that $\pi_*\mathcal O_X(\ulcorner 
\mathbf A^{\ast}(X, B)\urcorner +jb \overline {H})\subseteq
\pi_*\mathcal O_X(jbH)$ for every integer $j\geq j_0$, and 
\item[(3)] there is a positive integer $c$ such that 
$cH$ is Cartier and that $\mathcal O_T(cH):=\mathcal O_X(cH)|_T$ is 
$\pi$-generated, 
where 
$T=-\mathbf N(X, B)_X$.  
\end{itemize}
Then $H$ is $\pi$-semi-ample. 
\end{thm}
\begin{proof}
If $(X, B)$ is sub klt, then this follows from Theorem \ref{thm1}. 
By replacing $H$ by its multiple, 
we can assume that $b=1$, $j_0=1$, and $c=1$. 
Since $lH+\ulcorner \mathbf A^{\ast}_X\urcorner -T-(K_X+\{B\})
=lH-(K_X+B)$ is nef and big over $S$ for every positive 
integer $l$, 
we have the following commutative diagram by the Kawamata--Viehweg 
vanishing theorem: 
$$
\begin{CD}
\pi_*\mathcal O_X(lH+\ulcorner\mathbf A^{\ast}_X\urcorner)
@>>> \pi_*(\mathcal O_{T}(lH)\otimes\mathcal O_T(\ulcorner 
\mathbf A^{\ast}_X|_{T}\urcorner))@>>>0\\
@AA\text{$\cong$}A @AA\text{$\iota$}A\\
\pi_*\mathcal O_X(lH)
@>>{\alpha}>\pi_*\mathcal O_{T}(lH) . 
\end{CD}
$$ 
Thus, the natural inclusion $\iota$ is an isomorphism and $\alpha$ is 
surjective for every $l\geq 1$. In particular, 
$\pi_*\mathcal O_X(lH)\ne 0$ for every $l\geq 1$. 
The same arguments as in \cite[Proof of Theorem 3]{fukuda} 
show that $H$ is $\pi$-semi-ample. 
\end{proof}

The main purpose of this section is to prove Theorem \ref{abundantfukuda} 
below, which is a generalization of Theorem \ref{thm4} 
and Theorem \ref{bigfukuda}. 
The basic strategy of the proof is the same as that of 
Theorem \ref{thm4}. That is, by using Ambro's canonical 
bundle formula, we reduce it to the case when 
$H-(K_X+B)$ is nef and big. This is nothing but Theorem \ref{bigfukuda}. 

\begin{thm}\label{abundantfukuda}
Let $X$ be a non-singular variety and let $B$ be a $\mathbb Q$-divisor 
on $X$ such that $(X,B)$ is sub lc and 
$\Supp B$ is a simple normal crossing divisor. 
Let $\pi:X \to S$ be a proper 
morphism onto a variety $S$. Assume the following 
conditions{\em{:}} 
\begin{itemize}
\item[(a)] $H$ is a $\pi$-nef $\mathbb Q$-Cartier $\mathbb Q$-divisor on $X$, 
\item[(b)] $H-(K_X+B)$ is $\pi$-nef and $\pi$-abundant, 
\item[(c)] $\kappa (X_\eta, (aH-(K_X+B))_\eta)\geq 0$ 
and $\nu(X_\eta, (aH-(K_X+B))_\eta)=\nu(X_\eta, (H-(K_X+B))_\eta)$ for 
some $a\in \mathbb Q$ with $a>1$, 
where $\eta$ is the generic point of $S$,  
\item[(d)] let $f:Y\to Z$ be the Iitaka 
fibration with respect to $H-(K_X+B)$ over
$S$. We assume that there 
exists a proper birational morphism $\mu: Y\to X$
and put $K_Y+B_Y=\mu ^*(K_X+B)$. In this setting, we assume 
$\rank f_*\mathcal O_Y(\ulcorner \mathbf A^{\ast}
(Y, B_Y)\urcorner)=1$, 
\item[(e)] $($Saturation condition$)$
there exist positive integers $b$ and $j_0$ such that
$bH$ is Cartier and $\pi_*\mathcal O_X(\ulcorner
\mathbf A^{\ast}(X, B)\urcorner+jb\overline H)\subseteq \pi_*
\mathcal O_X(jbH)$ for every positive integer $j\geq j_0$, and 
\item[(f)] there is a positive integer $c$ such that 
$cH$ is Cartier and that 
$\mathcal O_T(cH):=\mathcal O_X(cH)|_T$ is 
$\pi$-generated, 
where 
$T=-\mathbf N(X, B)_X$.
\end{itemize}
Then $H$ is $\pi$-semi-ample. 
\end{thm}
\begin{proof}
If $H-(K_X+B)$ is big, then this follows from Theorem \ref{bigfukuda}. 
So, we can assume that $H-(K_X+B)$ is not big. 
Form now on, we use the notation 
in the proof of Theorem \ref{thm2} in \cite[Section 2]{fujino-kawamata}. 
We just explain how to modify that proof. 
Let us recall the commutative diagram 
$$
\begin{CD}
Y@>{f}>>Z\\
@V{\mu}VV @VV{\varphi}V\\
X@>>{\pi}>S
\end{CD}
$$
in the proof of \cite[Theorem 1.1]{fujino-kawamata}. 
For the details, see \cite[Section 2]{fujino-kawamata}. 
We start with the following obvious lemma. 
\begin{lem}\label{6.3}
We put $T'=-\mathbf N(X, B)_Y$. 
Then $\mu (T')\subset T$. 
Therefore, $\mathcal O_{T'}(cH_Y):=\mathcal O_Y(cH_Y)|_{T'}$ 
is $p$-generated, where 
$p=\pi\circ \mu$. 
\end{lem}
\begin{lem}
If $f(T')=Z$, then $H_Y$ is $p$-semi-ample. 
In particular, $H$ is $\pi$-semi-ample. 
\end{lem}
\begin{proof}
There exists an irreducible component $T'_0$ of $T'$ such 
that $f(T'_0)=Z$. Since $(H_Y)|_{T'_0}\sim _{\mathbb Q} (f^*D)|_{T'_0}$ 
is $p$-semi-ample, 
$D$ is $\varphi$-semi-ample. 
This implies that $H_Y$ is $p$-semi-ample and 
$H$ is $\pi$-semi-ample. 
\end{proof}
Therefore, we can assume that 
$T'$ is not dominant onto $Z$. 
Thus $\mathbf A(Y, B_Y)=\mathbf A^{\ast} (Y, B_Y)$ over the 
generic point of $Z$. 
Equivalently, $(Y, B_Y)$ is sub klt over the 
generic point of $Z$. 
By applying Ambro's result as in the proof 
of Theorem \ref{thm2} in \cite[Section 2]{fujino-kawamata}, we have the 
properties (1), (3), (5), 
\begin{itemize}
\item[(2$'$)] $(Z, B_Z)$ is sub lc, 
\item[(4$'$)] $\varphi_*\mathcal O_Z(\ulcorner 
\mathbf A^{\ast}(Z, B_Z)\urcorner +j\overline D)
\subseteq \varphi _*\mathcal O_Z(jD)$ for every 
positive integer 
$j$, 
\item[(6)] $Y$ and $Z$ are non-singular and $\Supp 
B_Y$ and $\Supp B_Z$ are simple normal 
crossing divisors, and 
\item[(7)] $\mathcal O_{T''}(D):=\mathcal O_Z(D)|_{T''}$ is 
$\varphi$-semi-ample where $T''=-\mathbf N(Z, B_Z)_Z$.  
\end{itemize}
Once the above conditions were satisfied, $D$ is $\varphi$-semi-ample 
by Theorem \ref{bigfukuda}. Therefore, $H$ is $\pi$-semi-ample. 
So, all we have to do is to check the 
above conditions. 
The conditions (1), (2$'$), (3), (5), (6) are satisfied 
by Ambro's result. 
We note that $\rank f_*\mathcal O_Y(\ulcorner \mathbf A(Y, B_Y)\urcorner) 
=\rank f_*\mathcal O_Y(\ulcorner \mathbf A^{\ast}(Y, B_Y)\urcorner)=1$. 
By the same computation as in \cite[Lemma 9.2.2 and 
Proposition 9.2.3]{ambro3}, 
we have the following 
lemma. 
\begin{lem}
$\mathcal O_Z(\ulcorner \mathbf A^{\ast} (Z, B_Z)\urcorner+j\overline 
D)\subseteq f_*\mathcal O_Y(\ulcorner 
\mathbf A^{\ast}(Y, B_Y)\urcorner +j\overline {H_Y})$ for 
every integer $j$. 
\end{lem}
Thus, we have (4$'$) by the saturation condition (e) (for the 
details, see the proof of Theorem \ref{thm2} 
in \cite[Section 2]{fujino-kawamata}, and Lemma \ref{lem321}). 
By definition, we have 
$lH_Y+\ulcorner \mathbf A^{\ast}_Y\urcorner -T'-(K_Y+\{B_Y\}) 
\sim_{\mathbb Q}f^*((l-1)D+M_0)$. 
Note that $(l-1)D+M_0$ is $\varphi$-nef and 
$\varphi$-big for $l\geq 1$. 
By the Koll\'ar type injectivity theorem, 
$$
p_*\mathcal O_Y(lH_Y+\ulcorner \mathbf A^{\ast}_Y\urcorner -T')
\to p_*\mathcal O_Y(lH_Y+\ulcorner \mathbf A^{\ast}_Y\urcorner)
$$ 
is injective for $l\geq 1$. 
Note that the above injectivity 
can be checked easily by \cite[Theorem 1.1]{fujino-trans}. 
Here, we used the fact that $f(T')\subsetneq Z$. So, we have the 
following commutative diagram: 
$$
\begin{CD}
p_*\mathcal O_Y(lH_Y+\ulcorner\mathbf A^{\ast}_Y\urcorner)
@>>> p_*(\mathcal O_{T'}(lH_Y)\otimes\mathcal O_{T'}(\ulcorner 
\mathbf A^{\ast}_Y|_{T'}\urcorner))@>>>0\\
@AA\text{$\cong$}A @AA\text{$\iota$}A\\
p_*\mathcal O_Y(lH_Y)
@>>{\alpha}>p_*\mathcal O_{T'}(lH_Y) . 
\end{CD}
$$ 
The isomorphism of the left vertical arrow 
follows from the saturation condition (e). 
Thus, the natural inclusion $\iota$ is an isomorphism 
and $\alpha$ is surjective for $l\geq 1$. 
In particular, the relative base locus of $lH_Y$ is disjoint from 
$T'$ if $\mathcal O_{T'}(lH_Y)$ is $p$-generated. 
On the other hand, $H_Y\sim _{\mathbb Q}f^*D$. 
Therefore, by Lemma \ref{6.3}, $\mathcal O_{T''}(D)$ is $\varphi$-semi-ample since 
$T''\subset f(T')$. 
So, we obtain the condition (7). We complete the 
proof of Theorem \ref{abundantfukuda}. 
\end{proof}

As a corollary of 
Theorem \ref{abundantfukuda}, we obtain a slight generalization of 
Fukuda's result (cf.~\cite[Proposition 3.3]{fukuda2}). 
Before we explain the corollary, let us recall 
the definition of {\em{non-klt loci}}. 

\begin{defn}[Non-klt locus] 
Let $(X, B)$ be an lc pair. 
We consider the closed subset 
$$
\Nklt(X, B)=\{x\in X \,|\, (X, B)\  \text{is not klt at}\  x\}
$$ 
of $X$. 
We call $\Nklt(X,B)$ the {\em{non-klt locus}} of $(X, B)$. 
\end{defn}

\begin{cor}\label{fuku-main}
Let $(X, B)$ be an lc pair and 
let $\pi:X \to S$ be a proper 
morphism onto a variety $S$. Assume the following 
conditions{\em{:}} 
\begin{itemize}
\item[(a)] $H$ is a $\pi$-nef $\mathbb Q$-Cartier $\mathbb Q$-divisor on $X$, 
\item[(b)] $H-(K_X+B)$ is $\pi$-nef and $\pi$-abundant, 
\item[(c)] $\kappa (X_\eta, (aH-(K_X+B))_\eta)\geq 0$ 
and $\nu(X_\eta, (aH-(K_X+B))_\eta)=\nu(X_\eta, (H-(K_X+B))_\eta)$ for 
some $a\in \mathbb Q$ with $a>1$, 
where $\eta$ is the generic point of $S$,  
\item[(f)] there is a positive integer $c$ such that 
$cH$ is Cartier and that $\mathcal O_T(cH):=\mathcal O_X(cH)|_T$ is 
$\pi$-generated, 
where $T=\Nklt(X, B)$ is the non-klt locus of $(X, B)$.  
\end{itemize}
Then $H$ is $\pi$-semi-ample. 
\end{cor} 
The readers can find applications of this corollary in \cite{fukuda2} 
and \cite{fujino-finite}. 
\begin{proof}
Let $h:X'\to X$ be a resolution such that 
$\Exc(h)\cup \Supp h^{-1}_*B$ 
is a simple normal crossing divisor and 
$K_{X'}+B_{X'}=h^*(K_X+B)$. 
Then $H_{X'}=h^*H$, $(X', B_{X'})$, 
and $\pi'=\pi\circ h:X'\to S$ satisfy the assumptions (a), 
(b), and (c) in Theorem \ref{abundantfukuda}. By the 
same argument as in the proof of 
\cite[Lemma 2.3]{fujino-kawamata}, we obtain 
$\rank f_*\mathcal O_Y(\ulcorner \mathbf A^{\ast}(Y, B_Y)\urcorner)=1$, 
where $f:Y\to Z$ is the Iitaka fibration as in (d) in Theorem 
\ref{abundantfukuda}. 
Note that $\ulcorner \mathbf A^{\ast}(Y, B_Y)\urcorner$ is 
effective and exceptional over $X$. 
Since $B$ is effective, 
$\ulcorner \mathbf A^{\ast}(X, B)\urcorner$ is effective 
and exceptional over 
$X$, 
$$\pi'_*\mathcal O_{X'}(\ulcorner 
\mathbf A^{\ast}(X', B_{X'})\urcorner+jb\overline {H_{X'}})
\subseteq \pi'_*\mathcal O_{X'}(jbH_{X'})$$ for 
every integer $j$, 
where $b$ is a positive integer such that 
$bH$ is Cartier. 
So, the saturation condition (e) in Theorem 
\ref{abundantfukuda} is satisfied. 
Finally, $\mathcal O_{T'}(cH_{X'}):=\mathcal O_{X'}(cH_{X'})|_{T'}$ is 
$\pi'$-generated, 
where $T'=-\mathbf N(X, B)_{X'}$, 
by the assumption (f) and the fact that 
$h(T')\subset T$. 
So, the condition (f) in Theorem \ref{abundantfukuda} 
for $H_{X'}$ and $(X', B_{X'})$ is satisfied. 
Therefore, $H_{X'}$ is $\pi'$-semi-ample by Theorem 
\ref{abundantfukuda}. Of course, 
$H$ is $\pi$-semi-ample.  
\end{proof}

\begin{rem}
(i) It is obvious that 
$\Supp (-\mathbf N(X, B)_X)\subseteq \Nklt (X, B)$. 
In general, $\Supp (-\mathbf N(X, B)_X)\subsetneq 
\Nklt (X, B)$. In particular, $\Nklt (X, B)$ is not necessarily 
of pure codimension one in $X$. 

(ii) If $(X, B)$ is dlt, then $\Nklt (X, B)=\Supp (-\mathbf N(X, B)_X)
=\llcorner B\lrcorner$. 
Therefore, if $(X, B)$ is dlt and 
$S$ is a point, then Corollary \ref{fuku-main} is 
nothing but Fukuda's result \cite[Proposition 3.3]{fukuda2}. 
\end{rem}
 
By combining Corollary \ref{fuku-main} with 
\cite[Theorem 1.5]{gongyo}, we obtain the 
following result. 

\begin{cor}
Let $(X, B)$ be a projective dlt pair 
such that 
$\nu(K_X+B)=\kappa (K_X+B)$ and 
that $(K_X+B)|_{\llcorner B\lrcorner}$ is numerically trivial. 
Then $K_X+B$ is semi-ample. 
\end{cor}

We close this section with a remark. 

\begin{rem}
We can easily generalize Theorem \ref{abundantfukuda} 
and Corollary \ref{fuku-main} to 
varieties in class $\mathcal C$ by suitable modifications. 
We omit details here. See \cite[Section 4]{fujino-kawamata}. 
\end{rem}

\section{Base point free theorems for pseudo-klt pairs}\label{sec-exc}

In this section, we generalize the Kawamata--Shokurov base point 
free theorem and Kawamata's theorem:~Theorem \ref{thm2} for {\em{klt 
pairs}} to {\em{pseudo-klt pairs}}. 
We think that 
our formulation is useful when we 
study lc centers (see Proposition \ref{7878}). 
First, we introduce the notion of 
{\em{pseudo-klt pairs}}. 

\begin{defn}[Pseudo-klt pair]\label{quasiklt}
Let $W$ be a normal variety. Assume 
the following conditions: 
\begin{itemize}
\item[(1)] there exist a sub klt pair 
$(V, B)$ and a proper surjective morphism $f:V\to W$ 
with connected fibers, 
\item[(2)] $f_*\mathcal O_V(\ulcorner \mathbf A(V, B)\urcorner)\simeq 
\mathcal O_W$, and 
\item[(3)] there exists a $\mathbb Q$-Cartier $\mathbb Q$-divisor 
$\mathcal K$ on $W$ such that $K_V+B\sim _{\mathbb Q}f^*\mathcal K$. 
\end{itemize}
Then the pair $[W, \mathcal K]$ is called a {\em{pseudo-klt pair}}. 
\end{defn}

Although it is the first time that we use the name of 
{\em{pseudo-klt pair}}, the notion of pseudo-klt pair appeared 
in \cite{fuji2}, where we proved the cone and contraction 
theorem for {\em{pseudo-klt}} pairs (cf.~\cite[Section 4]{fuji2}). 
We note that all the fundamental theorems for the log minimal model program for 
pseudo-klt pairs can be proved 
by the theory of quasi-log varieties (cf.~\cite{ambroQ}, 
\cite{fujino-intro}, and \cite{fujino-book}). 

\begin{rem}
In Definition \ref{quasiklt}, we assume that 
$W$ is normal. However, 
the normality of $W$ follows from 
the condition (2) and the normality of $V$. 
Note that $\ulcorner \mathbf A(V, B)\urcorner$ is effective. 
\end{rem}
\begin{rem}
In the definition of pseudo-klt pairs, 
if $(V, B)$ is klt, then $f_*\mathcal O_V(\ulcorner 
\mathbf A(V, B)\urcorner)\simeq \mathcal O_W$ is 
automatically satisfied. 
It is because $\ulcorner \mathbf A(V, B)\urcorner$ is 
effective and exceptional over $V$. 
\end{rem}

We note that a pseudo-klt pair is a very 
special example of 
Ambro's quasi-log varieties (see \cite[Definition 4.1]{ambroQ}). 
More precisely, 
if $[V, \mathcal K]$ is a pseudo-klt pair, 
then we can easily check that $[V, \mathcal K]$ is a {\em{qlc pair}}. 
See, for example, \cite[Definition 3.1]{fujino-intro}. 
For the details of the theory of quasi-log varieties, see \cite{fujino-book}. 

\begin{thm}Let $[W, \mathcal K]$ be a 
pseudo-klt pair. 
Assume that $(V, B)$ is klt and $W$ is projective or that 
$W$ is affine. 
Then we can find an effective $\mathbb Q$-divisor 
$B_W$ on $W$ such that 
$(W, B_W)$ is klt and 
that 
$\mathcal K\sim _{\mathbb Q}K_W+B_W$. 
\end{thm}
\begin{proof}
When $(X, B)$ is klt and $W$ is projective, 
we can find $B_W$ by \cite[Theorem 4.1]{ambro-moduli}. 
When $W$ is affine, this theorem follows from \cite[Theorem 1.2]{fuji2}. 
\end{proof}

It is conjectured that 
we can always find an effective $\mathbb Q$-divisor 
$B_W$ on $W$ such that 
$(W, B_W)$ is klt and 
$\mathcal K\sim _{\mathbb Q} K_W+B_W$. 

\begin{say}[Examples] We collect basic examples 
of pseudo-klt pairs. 

\begin{ex}
A klt pair is a pseudo-klt pair. 
\end{ex}

\begin{ex}
Let $f:X\to W$ be a Mori fiber space. Then we 
can find a $\mathbb Q$-Cartier $\mathbb Q$-divisor 
$\mathcal K$ on $W$ such that 
$[W, \mathcal K]$ is a pseudo-klt pair. 
It is because we can find an effective $\mathbb Q$-divisor 
$B$ on $X$ such that $K_X+B\sim _{\mathbb Q, f}0$ and 
$(X, B)$ is klt. 
\end{ex}

\begin{prop}\label{7878}
An exceptional lc center $W$ of an 
lc pair $(X, B)$ is a pseudo-klt pair for some $\mathbb Q$-Cartier 
$\mathbb Q$-divisor $\mathcal K$ on $W$. 
\end{prop}
\begin{proof}
We take a resolution $g:Y\to X$ such that 
$\Exc (g)\cup g^{-1}_*B$ has a simple normal crossing 
support. 
We put $K_Y+B_Y=g^*(K_X+B)$. 
Then $-B_Y=\mathbf A(X, B)_Y=\mathbf A_Y=\mathbf A^{\ast}_Y+\mathbf 
N_Y$, where 
$\mathbf N_Y=-\sum _{i=0}^{k} E_i$. 
Without loss of 
generality, we can assume that 
$f(E)=W$ and $E=E_0$. 
By shrinking $X$ around $W$, we can assume 
that $\mathbf N_Y=-E$. 
Note that $R^1g_*\mathcal O_Y(\ulcorner \mathbf A^{\ast}_Y\urcorner-E)
=0$ by the Kawamata--Viehweg vanishing theorem since 
$\ulcorner \mathbf A^{\ast}_Y\urcorner-E=K_Y
+\{-\mathbf A^{\ast}_Y\}-g^*(K_X+B)$. 
Therefore, $g_*\mathcal O_Y(\ulcorner 
\mathbf A^{\ast}_Y\urcorner )\simeq \mathcal O_X\to 
g_*\mathcal O_E(\ulcorner \mathbf A^{\ast}_Y|_{E}\urcorner)$ is 
surjective. 
This implies that 
$g_*\mathcal O_E(\ulcorner \mathbf A^{\ast}_Y|_{E}\urcorner)\simeq 
\mathcal O_W$. 
In particular, $W$ is normal. 
If we put $K_E+B_E=(K_Y+B_Y)|_E$, 
then $(E, B_E)$ is sub klt and 
$\mathbf A^{\ast}_Y|_{E}=\mathbf A(E, B_E)_E=-B_E$. 
So, $g_*\mathcal O_E(\ulcorner \mathbf A(E, B_E)\urcorner)=
g_*\mathcal O_E(\ulcorner -B_E\urcorner)\simeq \mathcal O_W$. 
Since $K_E+B_E=(K_Y+B_Y)|_E$ and $K_Y+B_Y=g^*(K_X+B)$, 
we can find a $\mathbb Q$-Cartier 
$\mathbb Q$-divisor $\mathcal K$ on $W$ such that 
$K_E+B_E\sim _{\mathbb Q}g^*\mathcal K$. 
Therefore, $W$ is a pseudo-klt pair. 
\end{proof}

We give an important remark on {\em{minimal}} lc centers. 

\begin{rem}[Subadjunction for minimal lc center]
Let $(X, B)$ be a projective or affine lc pair and 
let $W$ be a minimal lc center 
of the pair $(X, B)$. 
Then we can find an effective $\mathbb Q$-divisor 
$B_W$ on $W$ such that 
$(W, B_W)$ is klt and 
$K_W+B_W\sim _{\mathbb Q}(K_X+B)|_W$. 
For the details, see \cite[Theorems 4.1, 7.1]{fujino-gongyo}. 
\end{rem}
\end{say}

The following theorem is the Kawamata--Shokurov base 
point free theorem for pseudo-klt pairs. 
We give a simple proof depending on Kawamata's 
positivity theorem. Although Theorem \ref{quasi2} seems to be 
contained in \cite[Theorem 7.2]{ambroQ}, 
there are no proofs of \cite[Theorem 7.2]{ambroQ} in 
\cite{ambroQ}. 

\begin{thm}\label{quasi2}
Let $[W, \mathcal K]$ be a pseudo-klt pair, let $\pi:W\to S$ be a proper 
morphism onto a variety $S$ and let $D$ be a $\pi$-nef Cartier 
divisor on $W$. 
Assume that $rD-\mathcal K$ is $\pi$-nef and $\pi$-big for 
some positive integer $r$. 
Then $mD$ is $\pi$-generated for every $m\gg 0$. 
\end{thm}
\begin{proof}
Without loss of generality, we can assume that 
$S$ is affine. 
By the usual technique (cf.~\cite[Theorem 1]{kawa2} and 
\cite[Theorem 1.2]{fuji2}), 
we have 
$$\mathcal K+\varepsilon (rD-\mathcal K)
\sim_{\mathbb Q}K_W+\Delta_W$$ such 
that $(W, \Delta_W)$ is klt for some 
sufficiently small rational number $0<\varepsilon \ll 1$ 
(see also \cite[Theorem 8.6.1]{kosurvey}). 
Then $rD-(K_W+\Delta_W)\sim _{\mathbb Q}
(1-\varepsilon) (rD-\mathcal K)$, which is $\pi$-nef 
and $\pi$-big. 
Therefore, $mD$ is $\pi$-generated for every $m\gg 0$ by the 
usual Kawamata--Shokurov base point free theorem. 
\end{proof}

The next theorem is the main theorem of this section. 
It is a generalization of Kawamata's theorem 
in \cite{kawamata} (cf.~Theorem \ref{thm2}) 
for pseudo-klt pairs. 

\begin{thm}\label{quasi1}
Let $[W, \mathcal K]$ be a pseudo-klt pair and let $\pi:W\to S$ be a proper 
morphism onto a variety $S$. 
Assume the following conditions{\em{:}} 
\begin{itemize}
\item[(i)] $H$ is a $\pi$-nef $\mathbb Q$-Cartier $\mathbb Q$-divisor on $W$, 
\item[(ii)] $H-\mathcal K$ is $\pi$-nef and 
$\pi$-abundant, and 
\item[(iii)]$\kappa (W_\eta, (aH-\mathcal K)_\eta)\geq 0$ 
and 
$\nu (W_\eta, (aH-\mathcal K)_\eta)=\nu 
(W_\eta, (H-\mathcal K)_\eta)$ 
for some $a\in \mathbb Q$ with 
$a>1$, where 
$\eta$ is the generic point of $S$. 
\end{itemize}
Then $H$ is $\pi$-semi-ample. 
\end{thm}
\begin{proof} 
By the definition, there 
exists a proper surjective morphism 
$f:V\to W$ from a sub klt pair $(V, B)$. 
Without loss of generality, we can assume that 
$V$ is non-singular 
and $\Supp B$ is a simple normal crossing 
divisor. 
By the definition, $f_*\mathcal O_V(\ulcorner -B\urcorner)
\simeq \mathcal O_W$. 
From now on, we assume that $H$ is Cartier by 
replacing it with 
its multiple. Then 
$f_*\mathcal O_V(\ulcorner -B\urcorner 
+jH_V)\simeq  \mathcal O_W(jH)$ by the projection formula 
for 
every integer $j$, where 
$H_V=f^*H$.  
By pushing it by $\pi$, we 
have 
\begin{align*}
p_*\mathcal O_V(\ulcorner 
\mathbf A(V, B)\urcorner+j\overline{H_V})
&=p_*\mathcal O_V(\ulcorner -B\urcorner+jH_V)\\&\simeq \pi_{*}
\mathcal O_W(jH)\\&\simeq p_*\mathcal O_V(jH_V)
\end{align*} 
for every integer $j$, where $p=\pi\circ f$. 
This is nothing but the saturation condition:~Assumption 
(e) in Theorem \ref{thm4}. 
We put $L=H-\mathcal K$. 
We consider the Iitaka fibration with respect to 
$L$ over $S$ as in the 
proof of Theorem \ref{thm2} in \cite[Section 2]{fujino-kawamata}. 
Then we obtain the following commutative diagram: 
$$
\begin{CD}
V@= V\\
@V{f}VV @VVV\\
W@<\mu<< U\\ 
@V{\pi}VV @VV{g}V\\ 
S@<<{\varphi}<Z
\end{CD}
$$
where $g:U\to Z$ is the Iitaka fibration 
over $S$ and 
$\mu:U\to W$ is a birational 
morphism. 
Note that we can assume that 
$f:V\to W$ factors through 
$U$ by blowing up $V$. 
\begin{lem}\label{712}
$\rank h_*\mathcal O_V(\ulcorner \mathbf A(V, B)\urcorner)=1$, 
where $h:V\to U\to Z$. 
\end{lem}
\begin{proof}[Proof of {\em{Lemma \ref{712}}}] 
This proof is essentially the same as that of 
\cite[Lemma 2.3]{fujino-kawamata}. 
First, we can assume that 
$S$ is affine. 
Let $A$ be an ample divisor on $Z$ such 
that $h_*\mathcal O_V(\ulcorner \mathbf A(V, B)\urcorner)\otimes 
\mathcal O_Z(A)$ is $\varphi$-generated. 
We note that we can assume that 
$\mu^*L\sim_{\mathbb Q}g^*M$ since 
$L$ is $\pi$-nef and $\pi$-abundant, 
where $M$ is a $\varphi$-nef and $\varphi$-big 
$\mathbb Q$-divisor on $Z$. 
If we choose a large and divisible 
integer $m$, 
then $\mathcal O_Z(A)\subset \mathcal O_Z(mM)$. 
Thus 
\begin{eqnarray*}
&&\varphi_*(h_*\mathcal O_V(\ulcorner \mathbf A(V,B)\urcorner)
\otimes \mathcal O_Z(A)) 
\\&\subseteq &
\varphi_*(h_*\mathcal O_V(\ulcorner \mathbf A(V, B)\urcorner)
\otimes \mathcal O_Z(mM))
\\
&\simeq & p_*\mathcal O_V(\ulcorner 
\mathbf A(V, B)\urcorner+m\overline{f^*L})
\\&\simeq & \pi_{*}\mathcal O_W(mL)
\\&\simeq & \varphi_*\mathcal O_Z
(mM). 
\end{eqnarray*}
Therefore, we have 
$\rank h_*\mathcal O_V(\ulcorner \mathbf 
A(V, B)\urcorner)\leq 1$.  
Since $\mathcal O_Z\subset h_*\mathcal O_V\subset 
h_*\mathcal O_V(\ulcorner \mathbf A(V, B)\urcorner)$, 
we obtain $\rank h_*\mathcal O_V(\ulcorner \mathbf 
A(V, B)\urcorner)=1$ 
\end{proof}
Note that $h:V\to Z$ is the Iitaka fibration with respect to $f^*L$ over 
$S$. The assumption (c) in Theorem \ref{thm4} easily follows 
from (iii). Thus, by Theorem \ref{thm4}, we have 
that $H_V$ is $p$-semi-ample. 
Equivalently, $H$ is $\pi$-semi-ample. 
\end{proof}

The final theorem of this paper is a base point free theorem for 
minimal lc centers. 

\begin{thm}\label{final-thm}
Let $(X, B)$ be an lc pair and let $W$ be a minimal 
lc center of $(X, B)$. 
Let $\pi:W\to S$ be a proper morphism 
onto a variety $S$. 
Assume the following conditions{\em{:}} 
\begin{itemize}
\item[(i)] $H$ is a $\pi$-nef $\mathbb Q$-Cartier 
$\mathbb Q$-divisor on $W$, 
\item[(ii)] $H-(K_X+B)|_W$ is $\pi$-nef and 
$\pi$-abundant, and 
\item[(iii)] $\kappa (W_\eta, 
(aH-(K_X+B))|_{W_\eta})\geq 0$ and 
$\nu (W_\eta, (aH-(K_X+B))|_{W_\eta})
=\nu (W_\eta, (H-(K_X+B))|_{W_\eta})$ 
for some $a\in \mathbb Q$ with 
$a>1$, where 
$\eta$ is the generic point of $S$. 
\end{itemize}
Then $H$ is $\pi$-semi-ample. 
\end{thm}
\begin{proof}
Let $f:Y\to X$ be a dlt blow-up such that 
$K_Y+B_Y=f^*(K_X+B)$ 
(see, for example, \cite[Theorem 10.4]{fujino-fundamental}). 
Then we can take a minimal lc center $Z$ of $(Y, B_Y)$ such that 
$f(Z)=W$. Note that $K_Z+B_Z=(K_Y+B_Y)|_Z$ is klt. 
We also note that $W$ is normal (see, for example, \cite[Theorem 2.4 (4)]{fujino-non} 
or \cite[Theorem 9.1 (4)]
{fujino-fundamental}).  
Let 
$$
\begin{CD}
f:Z@>{g}>> V @>{h}>> W
\end{CD} 
$$ 
be the Stein factorization of $f:Z\to W$. 
Then $[V, h^*((K_X+B)|_W)]$ is a pseudo-klt pair by 
$g:(Z, B_Z)\to V$. We note that 
$H$ is $\pi$-semi-ample if and only if $h^*H$ is $\pi\circ h$-semi-ample. 
By Theorem \ref{quasi1}, $h^*H$ is semi-ample over $S$. 
We finish the proof. 
\end{proof}
%%%%%%%%%%%%%%%%%%%%%%%%%%%%%%%%%
\ifx\undefined\bysame
\newcommand{\bysame|{leavemode\hbox to3em{\hrulefill}\,}
\fi

\end{document}